\documentclass[11pt]{article}

\usepackage{mathtools,theorem,amssymb,amsmath}

\usepackage{graphicx}

\usepackage{hyperref}

\newtheorem{theorem}{Theorem}[section]

\theorembodyfont{\rm}

\newtheorem{remark}[theorem]{Remark}

\def\doubleprod#1#2{\ooalign{$#1\prod$\cr$#1\coprod$\cr}}
\DeclareMathOperator*{\rprod}{\mathpalette\doubleprod\relax}

\advance \topmargin by -\headheight
\advance \topmargin by -\headsep
\def \proof{\bigbreak\noindent{\it Proof.\ \ }}

\author{J.-P. Allouche \\
CNRS, IMJ-PRG \\
Sorbonne Universit\'e \\
4 Place Jussieu \\
F-75252 Paris Cedex 05 (France)
}

\title{The zeta-regularized product of odious numbers}

\date{ }

\begin{document}

\maketitle

\begin{abstract}
What is the product of all {\em odious} integers, i.e., of all integers whose binary expansion
contains an odd number of $1$'s? Or more precisely, how to define a product of these integers
which is not infinite, but still has a ``reasonable'' definition? We will answer this question by
proving that this product is equal to $\pi^{1/4} \sqrt{2 \varphi e^{-\gamma}}$, where $\gamma$  and
$\varphi$ are respectively the Euler-Mascheroni and the Flajolet-Martin constants.
\end{abstract}

% zeta-regularized product; odious numbers; zeta functions; Thue-Morse sequence

\begin{center}
{\it \hfill --- Dedicated to Joseph Kung} 
\end{center}

\section{Introduction}

Extending or generalizing ``simple'' notions is a basic activity in mathematics. This involves
trying to give some sense to an {\em a priori} meaningless formula, like $\sqrt{-1}$, $1/0$,
$\sum_{n \geq 1} n$, etc. Among these attempts is the question of ``assigning a reasonable value'' 
to an infinite product of increasing positive real numbers. This question arises for example when
trying to define ``determinants'' for operators on infinite-dimensional vector spaces.
One possible approach is to define ``zeta-regularization'' (see the definition in Section~\ref{def-ex} 
below). The literature on the subject is vast, going from theoretical aspects to explicit computations
in mathematics but also in physics (see, e.g., \cite{hawking, SZ}): we will ---of course--- not give a 
complete view of the existing references, but rather restrict to a few ones to allude to general
contexts where these infinite products take place. Our purpose is modest: to give the value
of an infinite arithmetic product (namely the product of all {\em odious} integers, i.e., of all those 
integers whose binary expansion contains an odd number of $1$'s, see, e.g., \cite{oeis})
$$
1 \times 2 \times 4 \times 7 \times 8 \times 11 \times 13 \times 14 \times 16 \times 19 
\times 21 \times 22 \times 25 \times 26 \times 28 \times 31 \times 32 \times \ldots 
$$ 
% and to show how this can be extended to 
% sets of integers having an automatic characteristic function.

\section{Definitions. First properties. Examples}\label{def-ex}

\subsection{Definitions}

The remark that 
$$
\log \prod_{i=1}^n \lambda_i = 
- \left(\frac{\rm d}{\rm ds} \sum_{i=1}^n \frac{1}{(\lambda_i)^s}\right)_{s=0}
$$
suggests a way of defining the infinite product of a sequence $(\lambda_i)_{i \geq 1}$ of positive
numbers (see, e.g., \cite{QHS, SABK}) by means of {\em zeta-regularization}: suppose that the 
Dirichlet series $\zeta_{\Lambda}(s) := \sum_i 1/\lambda_i^s$ converges when the real part of $s$ is 
large enough, that it has a meromorphic continuation to the whole complex plane, and that it has no 
pole at $0$, then the {\em zeta-regularized} product of the $\lambda_i$'s is defined by
$$
{\rprod_{i=1}^{\infty}}\lambda_i := e^{-\zeta_{\Lambda}'(0)}
$$ 
(this definition clearly coincides with the usual product when the sequence of $\lambda_i$'s is finite).
If $\zeta_{\Lambda}$ has as pole at $0$, there is a slight generalization of the definition above
(see \cite{illies, KW2}, also see \cite{KKSW}):
$$
{\rprod_{i=1}^{\infty}}\lambda_i := e^{-\underset{s=0}{\rm Res}\frac{\zeta_{\Lambda}(s)}{s^2}}
$$
where $\underset{s=0}{\rm Res} \ g(s)$ stands for the residue at $0$ of the function $g$.  
To put some general deep context about this definition, in particular about infinite determinants,
the reader can consult \cite{deninger, manin, voros}.

\subsection{First properties}\label{properties}

The following equalities hold

\begin{itemize}

\item[ ]{$*$} For all $N \geq 1$, \
$\displaystyle\rprod_{i=1}^{\infty}\lambda_i = 
\prod_{i=1}^N \lambda_i \ \left(\rprod_{i=N+1}^{\infty} \lambda_i\right)$.

\item[ ]{$*$}  For $a > 0$, \  $\displaystyle\rprod_{i=1}^{\infty} (a \lambda_i) = 
\left(\rprod_{i=1}^{\infty}\lambda_i \right) a^{\zeta_{\lambda}(0)}$

\item[ ]{$*$} If $A$ and $B$ form a partition of the positive integers, then
$\displaystyle\rprod_{i=1}^{\infty}\lambda_i = \displaystyle\rprod_{i \in A} \lambda_i
\displaystyle\rprod_{i \in B}\lambda_i$.

\end{itemize}

\subsection{Examples}

One can find several examples in the literature, (taken, e.g., from 
\cite{KKSW, KW1, KW2, QHS, yoshimoto} or deduced from properties in Section~\ref{properties}):

$$
\begin{array}{ll}
&\displaystyle\rprod_{n \geq 0}(n+x) = \frac{\sqrt{2\pi}}{\Gamma(x)} 
\ \ \ \mbox{\rm (for $x > 0$; Lerch formula $[*]$)} \\ 
&\displaystyle\rprod_{n \geq 1} n =  \sqrt{2 \pi} \ \ \mbox{\rm (Lerch formula for $x = 1$)}, \ 
\mbox{\rm hence, e.g.,} \ \displaystyle\rprod_{n \geq 1} (2n) = \sqrt{\pi}, \ \mbox{\rm and} \ 
\displaystyle\rprod_{n \geq 1} (2n+1) = \sqrt{2} \\
&\displaystyle\rprod_{n \geq 0} (n^2+1) = e^{\pi} - e^{-\pi} 
\ \ \mbox{\rm (general Lerch formula $[*]$)} \\
&\displaystyle\rprod_{n \geq 0} (n^2-n+1) = e^{\pi \frac{\sqrt{3}}{2}} + e^{-\pi \frac{\sqrt{3}}{2}} 
\ \ \mbox{\rm $[*]$} \\
&\displaystyle\rprod_{n \geq 0} (n^4+1) = 2 \left(\cosh(\pi\sqrt{2}) - \cos(\pi\sqrt{2})\right) 
\ \ \mbox{\rm $[*]$} \\
&\displaystyle\rprod_{n \geq 1} n^n = e^{-\zeta'(-1)} = A e^{-1/12}  \ \ \ 
(\mbox{\rm $A =  1.2824271 \ldots$ is the Glaisher-Kinkelin constant $[**]$)} \\
&\displaystyle\rprod_{n \geq 0} a^n = a^{-1/12} \ \ \mbox{\rm (for $a > 1$, $[***]$)} \\
&\!\displaystyle\rprod_{n \in Sq}  n = 2 \pi  \ \ \mbox{\rm (where $Sq$ is the set of squarefree 
positive integers; compare with $\displaystyle\rprod_{n \geq 1} n^2 = 2 \pi$)}\\
\end{array}
$$

\noindent
$[*]$ The original formula proved by Lerch (cited in \cite[p.~941--942]{KW1}) reads
$$
\frac{\rm d}{{\rm d}s} \left(\sum_{n \geq 0} \frac{1}{(n+x)^s}\right)_{s=0} =  \ \
\log \frac{\Gamma(x)}{\sqrt{2\pi}}\cdot
$$
Actually Lerch proved (cited in \cite[Equality (2), p.~942]{KW1}) the general formula
$$
\rprod_{n \geq 0} ((n+x)^2 + y^2) = \frac{2\pi}{\Gamma(x+iy)\Gamma(x-iy)}
$$
which of course implies the classical Lerch formula. This general Lerch formula was generalized
in \cite{KW1} where $((n+x)^2 + y^2)$ is replaced with $((n+x)^m + y^m)$.

\medskip

\noindent
$[**]$ Recall that the Glaisher-Kinkelin constant $A =  1.2824271 \ldots$ can be defined in several 
ways (see, e.g., \cite{vgorder} and the references therein):
$$
A = \lim_{n \to \infty} 
\frac{1^1 \times 2^2 \times 3^3 \times... \times n^n}{n^{n^2 + n/2 + 1/12}} \ e^{n^2/4}
= e^{\frac{1}{12} - \zeta'(-1)} = (2\pi)^{\frac{1}{12}} e^{\frac{\gamma}{12} - \frac{\zeta'(2)}{2\pi^2}} 
= \left(2\pi e^{\gamma} \prod_{p \ {\rm prime}} p^{\frac{1}{p^2 -1}}\right)^{1/12}.
$$

\medskip

\noindent
[***] As indicated in \cite{KW2} one may recall that formally
$\displaystyle\sum_{n \geq 1} n = \zeta(-1) = -1/12$. 

\bigskip

Another example of zeta-regularized product is given by the Fibonacci numbers $(F_n)_{n \geq 1}$ 
in \cite{kitson} where $\rprod_{n =1}^{\infty} F_n$  is computed in terms of the Fibonacci factorial 
constant and the golden ratio or in terms of the derivative of the Jacobi theta function of the first kind 
and the golden ratio. Of course the result needs the study of the Dirichlet series $\sum 1/F_n^s$: other 
similar Dirichlet series or zeta-regularized products are studied in \cite{BDR} and \cite{DSPR}.

\bigskip

Up to generalizing the notion of zeta-regularized product (``super-regularization''), one has
(\cite{MP}, also see \cite{smirnov}):
$$
\rprod_{p \ {\rm prime}} p = 4 \pi^2.
$$

\section{The zeta-regularized product of odious numbers}

In what follows, we let ${\mathcal O}$ denote the set of {\em odious} positive integers, i.e., of positive 
integers whose sum of binary digits is odd, and ${\mathcal E}$ the set of {\em evil} numbers,
i.e., of positive integers whose sum of binary digits is even. 
We let $a(n)$ denote the characteristic function of odious numbers (i.e., $a(n) = 1$ if the sum 
of binary digits of $n$ is odd, and $a(n)=0$ if it is even) and $\varepsilon_n = (-1)^{a(n)}$ (note 
that $\varepsilon_n = 1 - 2a(n)$). Clearly $\varepsilon_n$ is equal to $+1$ if the sum of the
binary digits of $n$ is even, and to $-1$ if the sum is odd; in other words 
$(\varepsilon_n)_{n \geq 1}$ is (up to its first term) the famous Thue-Morse sequence on 
the alphabet $\{-1, +1\}$ (see, e.g., \cite{AS}).

\begin{theorem}\label{prod-odious}
We have with the notation above, with 
$Q = \displaystyle\prod_{n \geq 1} \left(\frac{2n}{2n+1}\right)^{\varepsilon _n}$,
$$
\begin{array}{ll}
&\displaystyle\rprod_{n \in {\mathcal O}} = (2\pi)^{1/4} Q^{-1/2} \\
&\displaystyle\rprod_{n \in {\mathcal E}} = (2\pi)^{1/4} Q^{1/2}
\end{array}
$$
Also $Q = \displaystyle\frac{2^{-1/2} e^{\gamma}}{\varphi}$
where $\varphi$ is the Flajolet-Martin constant \cite{FM}.
\end{theorem}

\proof
For $\Re s > 1$ we have
$$
\zeta_{\mathcal O} (s) = \sum_{n \in {\mathcal O}} \frac{1}{n^s} 
= \sum_{n \geq 1} \frac{a(n)}{n^s} = \sum_{n \geq 1} \frac{1 - \varepsilon_n}{2n^s}
= \frac{1}{2} \zeta(s) - \frac{1}{2} g(s)
$$
where $\zeta$ is the 	Riemann zeta function and 
$g(s) := \displaystyle\sum_{n \geq 1} \frac{\varepsilon_n}{n^s}\cdot$ 

\noindent
But, by \cite[Theorem~1.2]{AC} with $q=2$ (also see \cite[Lemma~1]{FM}) 
$g$ can be analytically continued to the whole complex plane, and it satisfies, for all 
$s \in {\mathbb C}$,
\begin{equation}\label{dirichlet-morse}
g(s) = -1 + \sum_{k \geq 1} (-1)^{k+1} {s+k-1 \choose k} \frac{g(s+k)}{2^{s+k}}\cdot
\end{equation}
This implies that $g(0) = -1$ and 
$$
\frac{g(s)-g(0)}{s} = \sum_{k \geq 1} (-1)^{k+1} \frac{(s+1)(s+2)\ldots(s+k-1)}{k!} \frac{g(s+k)}{2^{s+k}}
$$
hence, by letting $s$ tend to $0$:
$$
g'(0) = \sum_{k \geq 1} (-1)^{k+1} \frac{g(k)}{k2^k}\cdot
$$
On the other hand, mimicking a computation in \cite[p.~534]{AC}, one has
$$
\sum_{k \geq 1} (-1)^{k+1} \frac{g(k)}{k2^k} = 
\sum_{k \geq 1} \frac{(-1)^{k+1}}{k2^k} \sum_{n \geq 1} \frac{\varepsilon_n}{n^k}
= \sum_{n \geq 1} \varepsilon_n \sum_{k \geq 1} \frac{(-1)^{k+1}}{k 2^k n^k}
= \sum_{n \geq 1} \varepsilon_n \log\left(1 + \frac{1}{2n}\right) = - \log Q
$$
where $\displaystyle Q := \prod_{n \geq 1} \left(\frac{2n}{2n+1}\right)^{\varepsilon_n}$.
So that $g'(0) = - \log Q$. We thus obtain
$$
\zeta_{\mathcal O}'(0) = \frac{1}{2} \zeta'(0) - \frac{1}{2} g'(0) = 
- \frac{1}{4} \log(2\pi) + \frac{1}{2} \log Q
$$ 
which finally yields
$$
\rprod_{n \in {\mathcal O}} n = e^{\frac{1}{4} \log(2\pi)  - \frac{1}{2} \log Q}
= (2\pi)^{1/4} Q^{-1/2}.
$$
Now
$$
\rprod_{n \in {\mathcal O}} n \rprod_{n \in {\mathcal E}} n = \rprod_{n \geq 1} n = \sqrt{2 \pi}
$$
which gives 
$$
\rprod_{n \in {\mathcal E}} n = (2\pi)^{1/4} Q^{1/2}.
$$
It remains to recall that the Flajolet-Martin constant $\varphi$ \cite{FM} is equal to
$$
\varphi := 2^{-1/2}e^{\gamma}\frac{2}{3} 
\prod_{n \geq 1} \left(\frac{(4n+1)(4n+2)}{4n(4n+3)}\right)^{\varepsilon_n} = 0.77351\ldots
$$
and thus (\cite[Section~6.8.1]{finch}, also see \cite{All-Bx}) 
$$
\varphi := \frac{2^{-1/2}e^{\gamma}}{Q} \  \ \mbox{\rm hence} \ \
Q = \frac{2^{-1/2} e^{\gamma}}{\varphi}\cdot
$$

\begin{remark}
Instead of considering the odious and evil numbers, one might have considered --in a rather
non-natural way-- the shifted odious and evil numbers, namely the sets
${\mathcal O_S} := \{n+1, \ n \in {\mathcal O }\}$ and ${\mathcal E_S} := \{n+1, \ n \in {\mathcal E }\}$.
Then $\rprod_{n \in {\mathcal O_S}} n = \exp(- \zeta_{\mathcal O_S}'(0))$.  
But, with the notation of the proof of Theorem~\ref{prod-odious}, and using that $\varepsilon_0 = 1$,
$$
\zeta_{\mathcal O_S}(s) = \sum_{n \in {\mathcal O}} \frac{1}{(n+1)^s} 
= \sum_{n \geq 1} \frac{a(n)}{(n+1)^s} = \frac{1}{2} \sum_{n \geq 1} \frac{1 - \varepsilon_n}{(n+1)^s}
= \frac{1}{2} \sum_{n \geq 0} \frac{1 - \varepsilon_n}{(n+1)^s} 
= \frac{1}{2} \zeta(s) - \frac{1}{2} f(s)
$$
\end{remark}
where $\displaystyle f(s) = \sum_{n \geq 0} \frac{\varepsilon_n}{(n+1)^s}$.
It was proved in \cite{AC} that this function $f$ can be analytically continued to
the whole complex plane, and that $f'(0) = \frac{\log 2}{2}$. Hence
$\zeta_{\mathcal O_S}'(0) = -\frac{1}{4} (\log(2 \pi) + \log 2) = -\frac{1}{4} \log 4\pi$.
This gives finally
$$
\rprod_{n \in {\mathcal O_S}} n = 2^{1/2} \pi^{1/4} \ \ \mbox{\rm and hence} \ \
\rprod_{n \in {\mathcal E_S}} n = \pi^{1/4}.
$$

\section{Beyond odious and evil}

The main tool used in the computation of the zeta-regularized product of odious or of evil
numbers above is the ``infinite functional equation''~(\ref{dirichlet-morse}). A multidimensional 
analog of this equation exists for ``automatic'' sequences \cite{AMFP}. We could expect
a similar result in the general case of a zeta-regularized product
of integers having an automatic characteristic function. What makes things more complicated
in the general case is that the involved zeta function occurs as a component of a vector of Dirichlet 
series satisfying an infinite functional equation, but that the zeta function itself does not necessarily 
satisfy such an equation. Furthermore this vector is meromorphic but not necessarily holomorphic 
on ${\mathbb C}$ (in particular $0$ might be a pole). As suggested by the referee, looking at the 
subcase given by the parity of block-counting sequences could well be the right generalization of 
the Thue-Morse case. Trying to follow this suggestion, we only arrived at a partial result for, e.g., 
the Golay-Shapiro (also called Rudin-Shapiro) sequence where, instead of considering the parity of 
the number of $1$'s in the binary expansions of integers for the Thue-Morse sequence, one considers 
the number of $11$'s in the binary expansions of integers. We hope to revisit these questions in the 
near future.

\bigskip

\noindent
{\bf Acknowledgments} We thank M. Marcus for having detected misprints in the first version of 
this paper. We also thank J.-Y. Yao, B. Saffari, and the referee for their constructive remarks.

\end{document}